\documentclass[11pt]{article}

\usepackage{amsmath}
\usepackage{amsfonts}
\usepackage{amsthm}
\usepackage{booktabs}
\usepackage{graphicx}
\usepackage{makeidx}
\usepackage{tocloft}
\usepackage{parskip}
\usepackage{verbatim}
\makeindex
\usepackage{listings}
\usepackage{url}
\usepackage[utf8x]{inputenc}
\usepackage[usenames]{color}

\definecolor{lg}{rgb}{0.9,0.9,0.9}
\definecolor{dg}{rgb}{0.3,0.3,0.3}

\title{Improving Non-linear Fits}

\author{Massimo Di Pierro \\
{\footnotesize School of Computing, DePaul University, 243 S Wabash Ave, Chicago IL 60604}}

\begin{document}

\maketitle

\begin{abstract}
In this notes we describe an algorithm for non-linear fitting which incorporates some of the features of linear least squares into a general minimum $\chi^2$ fit and provide a pure Python implementation of the algorithm. It consists of the variable projection method (varpro), combined with a Newton optimizer and stabilized using the steepest descent with an adaptative step. The algorithm includes a term to account for Bayesian priors. We performed tests of the algorithm using simulated data. This method is suitable, for example, for fitting with sums of exponentials as often needed in Lattice Quantum Chromodynamics. 
\end{abstract}

\section{Introduction}

Fitting is one the most common recurrent tasks in scientific computing. Fitting involves minimizing the $\chi^2$, i.e. weighted sum of squares of the difference between observed data and expected data~\cite{chi2}. The expected data is computing using a function $f$ which depends on some unknown parameters. The $\chi^2$ is minimized by varying these unknown parameters in order to find a minimum.

Usually there are two common approaches to fitting:
\begin{itemize}
\item If the function $f$ is linear in the unknown parameters, then the fitting problem has an exact solution which can be computed via the linear least squares.
\item If the function $f$ is nonlinear in its parameters, the minimization is carried out using a non-linear optimization algorithm based on the Newton method~\cite{newton}, the Steepest Descent~\cite{steepest}, or more sophisticated techniques such as the Levenberg-Marquadt.
\end{itemize}

There is also a {\it non-linear least squares} algorithm. It is different from the one described in this paper and it is equivalent to performing a Newton optimization where the first order derivatives of the fitting functions are computed analytically.

When the fitting function depends of both some linear parameters and some non-linear ones, the two methods above can be combined. This approach, described here, is known as VARiable PROjection method~\cite{varpro}. An implemention of it which uses the Levenberg-Marquadt algorithm for the non-linear part can be found in Netlib~\cite{netlib}. The approach consists of defining a function $g$ which takes as input the non-linear parameters of the original fitting function and returns the $\chi^2$. The function $g$, internally uses the least square algorithm to compute the linear coefficients exactly (as function of the non-linear coefficients and the input data) and uses this result to compute the $\chi^2$. In our implementaton the complete fit is performed by passing the function $g$ to a non-linear Newton optimizer stabilized using an adaptative steepest-descent. 

One application of this procedure is, for example, in lattice QCD, where typically one has to fit data points using a sum of exponentials of the form $a_i e^{-b_i t}$. $a_i$ and $b_i$ are the parameters to be determined. 

Fits like the one described become unstable when two of the $b$ coefficients are very close to each other. In order to avoid this situation we modified the algorithm to account for {\it a priori} constraints on the fitting parameters in the form of a Bayesian contribution to the $\chi^2$. This can be used to stabilize the fit~\cite{lepage}.

\section{Algorithm}

We want to fit data points of the form $(x_i,y_i \pm \delta y_i)$. 
We will do it by minimizing the $\chi^2$ defined as
\begin{equation}
\chi^2(\mathbf{a},\mathbf{b}) = \sum _i \left|
\frac{y_i - f(x_i,\mathbf{a},\mathbf{b})}{\delta y_i}
\right|^2
\end{equation}
The function $f$ is known but depends on unknown parameters $\mathbf{a}=(a_0,a_1,...)$ and $\mathbf{b}=(a_0,a_1,...)$. In terms of these parameters the function $f$ can be written as follows:
\begin{equation}
f(x,\mathbf{a},\mathbf{b}) = \sum_j a_j f_j(x, \mathbf{b})
\end{equation}
An example is the following one:
\begin{equation}
f(x,\mathbf{a},\mathbf{b}) = a_0 e^{-b_0 x} + a_1 e^{-b_1 x} + a_2 e^{-b_2 x} + ...
\end{equation}
The goal of our algorithm is to efficiently determine the parameters $\mathbf{a}$ and $\mathbf{b}$ which minimize the $\chi^2$.

We proceed by defining the following two quantities:

\begin{equation}
A(\mathbf{b}) = \left(
\begin{tabular}{cccc}
$f_0(x_0,\mathbf{b})/\delta y_0$ &
$f_1(x_0,\mathbf{b})/\delta y_0$ &
$f_2(x_0,\mathbf{b})/\delta y_0$ & ... \\
$f_0(x_1,\mathbf{b})/\delta y_1$ &
$f_1(x_1,\mathbf{b})/\delta y_1$ &
$f_2(x_1,\mathbf{b})/\delta y_1$ & ... \\
$f_0(x_2,\mathbf{b})/\delta y_2$ &
$f_1(x_2,\mathbf{b})/\delta y_2$ &
$f_2(x_2,\mathbf{b})/\delta y_2$ & ... \\
... & ... & ... & ...
\end{tabular}
\right)\qquad
\mathbf{z} = \left(
\begin{tabular}{c}
$y_0$ / $\delta y_0$ \\
$y_1$ / $\delta y_1$ \\
$y_2$ / $\delta y_2$ \\
...
\end{tabular}
\right)
\end{equation}

In terms of $A$ and $z$ the $\chi^2$ can be rewritten as

\begin{equation}
\chi^2(\mathbf{a},\mathbf{b}) = \left| A(\mathbf{b}) \mathbf{a} - \mathbf{z} \right|^2
\end{equation}

We can minimize this function (in $a$) using the linear least squares algorithm, exactly:

\begin{equation}
\mathbf{a}(\mathbf{b}) = (A(\mathbf{b})^t A(\mathbf{b}))^{-1}(A(\mathbf{b})^t \mathbf{z})
\end{equation}

We define a function which returns the minimum $\chi^2$ for a fixed input $\mathbf{b}$:

\begin{equation}
g(\mathbf{b}) =
\min_{\mathbf{a}}
\chi^2(\mathbf{a},\mathbf{b}) =
\chi^2(\mathbf{a}(\mathbf{b}),\mathbf{b}) = \left|
A(\mathbf{b})\mathbf{a}(\mathbf{b}) - \mathbf{z}
\right|^2 + \textrm{Bayesian}(\mathbf{b})
\end{equation}

Hence we have reduced the original problem to a simple problem by reducing the number of unknown parameters from $N_a+N_b$ to $N_b$.

We have arbitrarily added a term, {\tt Bayesian} to the definition of $\chi^2$. This term can be set to zero, but also, according to Lepage {\it et al.}~\cite{lepage}, we can choose it as follows:

\begin{equation}
\textrm{Bayesian}(\mathbf{b}) = \sum_i \left| \frac{b_i - \tilde b_i}{\delta b_i}\right|^2
\end{equation}
where $b_i \simeq \tilde b_i \pm \delta b_i$ is the {\it a priori} knowledge (or assumption) about the expected values for the $b$ parameters, according to Bayesian statistics. The Bayesian term originates from the fact that ordinary statistics concerns the computation of probability that data is compatible with a given model while here we are concerned with the opposite problem, the probability that a certain model (defined by the unknown parameters $a$ and $b$) is compatible with given data. In practice the effect of the {\tt Bayesian} term is that of stabilizing the fit by constraining the values within a range specified by a point and its uncertainty by increasing the $\chi^2$ quadratically when the fitting parameters depart from the priors.

The methodology described here (excluding the Bayesian contribution) is known as VARPRO~\cite{varpro} algorithm.

\section{Code}

The complete code is reported in the the Appendix. It contains the algorithm implemented in Python (use suggest version 2.7) and it uses the {\tt numpy} library for linear algebra. The code defines partial derivatives, the gradient, the Hessian, and the 1-norm of a matrix. These functions are necessary for the inner workings of the Newton algorithm. The code also includes a modified multi-dimensional Newton optimization algorithm that reverts to the steepest descent when the Newton step fails to reduce the value of the function passed as input. The step of the minimum descent is adjusted automatically: it starts with a guess and reduces the step until the $\chi^2$ is decreased by the move. This guarantes a decrease of the $\chi^2$ at each iteration. Although it never explodes, it may still fail. This can happen because the target precision is not reached within the maximum number of allowed iterations or because the Hessian becomes singular, or because $A$ matrix in the least squared fit becomes singular.

The main logic described in this notes is implemented in the function called {\tt fit}.

The input of the function {\tt fit} is:

\begin{itemize}
\item {\tt data} is a list of $(x,y,\delta y)$ points to be fitted.
\item {\tt fs} is a list of functions $f_i$ as described in the previous section. Each function must have the signature {\tt lambda b,x: ...} where $b$ is a list of $b$ coefficients.
\item {\tt b} is a list of $b$ values from which to start the optimization.
\item {\tt ap=1e-6} is the required target precision in the output.
\item {\tt rp=1e-4} is the required target relative precision.
\item {\tt ns=200} is the maximum number of steps performed by the Newton algorithm. If the Newton algorithm fails to converge within the required target absolute or relative precision within {\tt ns} steps, the algorithm raises an exception.
\item {\tt bayesian=None} can be set to a function of {\tt b} whose value is be added to the $\chi^2$ in the intermediate steps.
\end{itemize}

The function {\tt fit} returns a tuple containing $a$, $b$, the $\chi^2$ and the Hessian. $a$ and $b$ are lists.

\section{Examples}

\begin{figure}[t]
\begin{center}
\includegraphics[width=4in]{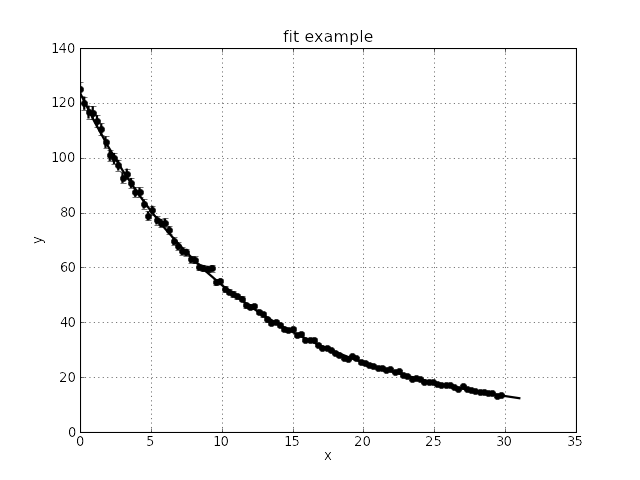}
\end{center}
\caption{Example of simulated data (100 points) and their fit with three exponentials\label{fig1}}
\end{figure}

\begin{figure}[t]
\begin{center}
\includegraphics[width=4in]{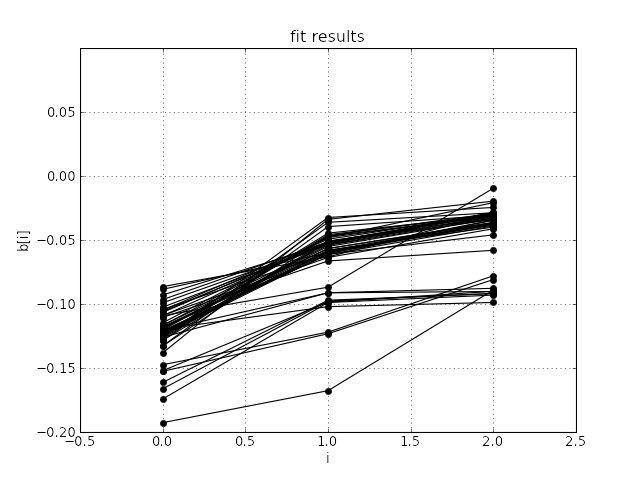}
\end{center}
\caption{Parallel coordinates plot for 50 simulated experiments. The data from each simulated experiment is fitted and the resulting $b_0$, $b_1$, $b_2$ are shown in the plot, connected by a line. The density of lines indicates the most likely range of results and their correlation. The results are compatible with the input used to generate the simulated experiments: $b_0=-0.10$, $b_1=-0.04$, $b_2=-0.02$.\label{fig2}}
\end{figure}

Here we report two examples of usage of the algorithm.
In the first one we consider a fitting function of the form
\begin{equation}
f(x) = a_0 x + a_1 x^2 + \frac{a_2}{x+b_0}
\end{equation}
We generated 100 points using the input values of $a_0=1$, $a_1=2$, $a_2=300$ and $b_0=10$ for $x=0.1,0.2,0.3,...10$. For each point we computed $y$ using the above formula and we added a random Gaussian noise with standard deviation equal to 1\% of $y$. We set  $\delta y = 1\%$ for each point. We then fitted the data using our method:

\begin{lstlisting}
>>> import random
>>> def noise(i,x,r=0.05): return (float(i),x*(1.0+random.gauss(0,r)),abs(x)*r)
>>> def f(x): return x+x**2+300.0/(x+10)
>>> data = [noise(0.1*i,f(0.1*i),0.01) for i in range(1,101)]
>>> fs = [(lambda b,x: x), 
          (lambda b,x: x*x), 
          (lambda b,x: 1.0/(x+b[0]))]
>>> a , b, chi2, H = fit(data,fs,b=[5])
>>> print a, b
[1.053.., 0.995..., 299.065...] [9.990...]
\end{lstlisting} 

The result shows that the fit converged within a few percent accuracy to the true values. Of course we do not expect the output to yield the exact same values which we used to generate the data, because we added the noise. The point here is that despite the noise, the fit is stable and consistent with expectation.

Notice that the notation {\tt lambda b,x: x*x} is nothing but the inline definition of an un-named function which takes {\tt b,x} and input and produces {\tt x*x} as output. Hence {\tt fs} is a list containing three functions.

In our second test we used a function which is the sum of three exponentials:
\begin{equation}
f(x) = a_0 e^{-b_0 x}+a_1 e^{-b_1 x}+a_2 e^{-b_2 x}
\end{equation}
We performed a series of 50 simulated experiments. In each experiment we generated 100 data points $x$ equally spaced between 0 and 32 and computed the corresponding $y$'s using the above formula with the input $a_0=100$, $a_1=10$, $a_2=1$ and $b_0 = -0.10$, $b_1 = -0.04$ and $b_2 = -0.02$ adding a 2\% Gaussian noise and assuming a 2\% error on each point. We then fit each simulated experiment. To improve the stability of the fit we added a Bayesian contribution to the $\chi^2$ which corresponds to priors $b_0 \simeq -0.11 \pm 0.04$, $b_1 \simeq -0.05 \pm 0.04$, and $b_2 \simeq -0.03 \pm 0.04$.

Here is the code for each simulated experiment:
\begin{lstlisting}
>>> from math import exp
>>> def f(x): return 100.0*exp(-0.10*x)+20.0*exp(-0.04*x)+4.0*exp(-0.02*x)
>>> xs = [0.3*i for i in range(100)]
>>> data = [noise(x,f(x),0.02) for x in xs]
>>> fs = [(lambda b,x: exp(b[0]*x)), 
          (lambda b,x: exp(b[1]*x)), 
          (lambda b,x: exp(b[2]*x))]
>>> def bayesian(b): return ((b[0]+0.11)**2+(b[1]+0.05)**2+(b[2]+0.03)**2)/0.04**2
>>> a, b, chi2, h = fit(data,fs,b=[-0.11,-0.05,-0.03],ns=100,bayesian=bayesian)
>>> print a, b
[98.041..., 10.795..., 16.281...] , [-0.103..., -0.052..., -0.029...]
\end{lstlisting}

The result of the fit is also shown in fig.~\ref{fig1}.

In this example the contribution of the third exponential is too small to be determined form data. We have ignored results that are more than two standard deviations off from priors and we have displayed the resulting fitting values for $b_0$, $b_1$ and $b_2$ in a parallel coordinates plot, fig.~\ref{fig2}. The X-axis represent the $i$ index of the $b_i$ parameters and the Y-axis the corresponding values. The lines connect the $b_0,b_1,b_2$ values for each simulated experiment. Higher density of lines represents the most likely values and their correlation.

We found that in about 10\% of our simulated experiments the algorithm failed to converge but this could be corrected by changing slightly the input starting point for the Newton optimizer.

\section{Conclusions}

In this notes we descibe the varpro method for non-linear fitting which uses linear least squares fitting to reduce the space of the parameters to be explored. We provided an implementation in Python which uses a stabilized Newton method the non-linear part of the fit and includes contribution from Bayesian priors. In the case of fitting with a sum of exponentials it reduces the size of the parameter space by a factor of two. We have provided working code and examples of fitting simulated data.

\subsection*{Acknowledgements}
We thank Carleton DeTar and Doug Toussaint for a useful discussion on this topic. We specially thank Martin Savage for pointing out to us that the method described here is known as the VARPRO algorithm, thus giving us the chance to properly acknowledge the original inventors~\cite{varpro}.

\newpage\section*{Appendix}

\begin{lstlisting}
### Developed by Massimo Di Pierro <mdipierro@cs.depaul.edu>
### License: BSD

from numpy import matrix
from numpy.linalg import *

def partial(f,i,h=1e-4):
    """
    definition of paritial derivative, df/dx_i
    """
    def df(x,f=f,i=i,h=h):
        x[i]+=h
        u = f(x)
        x[i]-=2*h
        v = f(x)
        x[i]+=h
        return (u-v)/2/h
    return df

def gradient(f, x, h=1e-4):
    """
    gradient of f in x
    """
    s = xrange(len(x))
    return matrix([[partial(f,r,h)(x)] for r in s])

def hessian(f, x, h=1e-4):
    """
    hessian of f in x
    """
    s = xrange(len(x))
    grad = [partial(f,r,h) for r in s]
    return matrix([[partial(grad[r],r,h)(x) for c in s] for r in s])

def norm(A):
    """
    defines norm of a matrix to check convergence
    """
    rows, cols = A.shape
    return max([sum(abs(A[r,c]) for r in xrange(rows)) \
                    for c in xrange(cols)])

def tolist(A):
    rows, cols = A.shape
    return [A[r,0] for r in xrange(rows)]

def optimize_newton_multi_imporved(f, x, ap=1e-6, rp=1e-4, ns=20):
    """
    Multidimensional Newton optimizer
    on failure is performs a steepest descent
    """
    fx = f(x)
    x = matrix([[element] for element in x])
    h = 10.0
    for k in xrange(ns):
        grad = gradient(f,tolist(x))
        (grad,H) = (gradient(f,tolist(x)), hessian(f,tolist(x)))
        if norm(H) < ap:
            raise ArithmeticError, 'unstable solution'
        (fx_old, x_old, x) = (fx, x, x-(1.0/H)*grad)
        fx = f(tolist(x))
        while fx>fx_old: # revert to steepest descent
            (fx, x) = (fx_old, x_old)
            n = norm(grad)
            (x_old, x) = (x, x - grad/n*h)
            (fx_old, fx) = (fx, f(tolist(x)))
            h = h/2
        h = norm(x-x_old)*2
        if k>2 and h/2<max(ap,norm(x)*rp):
            x = tolist(x)
            return x, hessian(f,x)
    raise ArithmeticError, 'no convergence'


def fit(data, fs, b, ap=1e-6, rp=1e-4, ns=200, bayesian=None):
    na = len(fs)
    def least_squares(b,data=data,fs=fs):
        A = matrix([[fs[k](b,x)/dy for k in xrange(na)] for (x,y,dy) in data])
        z = matrix([[y/dy] for (x,y,dy) in data])
        a = inv(A.T*A)*(A.T*z)
        chi2 = norm(A*a-z)**2
        return a, chi2
    def g(b,data=data,fs=fs,bayesian=bayesian):
        a, chi2 = least_squares(b, data, fs)
        if bayesian:
            chi2 += bayesian(b)
        return chi2
    b, H = optimize_newton_multi_imporved(g,b,ap,rp,ns)
    a, chi2 = least_squares(b,data,fs)
    return tolist(a), b, chi2, H
\end{lstlisting}

\end{document}